\documentclass [a4paper, 12pt] {amsart}

\usepackage {hyperref}

\usepackage [utf8] {inputenc}

\usepackage{mathtext}

\usepackage {version}

\usepackage {amsmath}
\usepackage {amssymb}
\usepackage {amsthm}

\usepackage {graphicx}

\usepackage {appendix}

\usepackage{mathtext}
\usepackage{amsmath}
\usepackage{amsfonts}
\usepackage{amssymb}
\usepackage{mathrsfs}
\usepackage{amsthm}

\usepackage {needspace}

\excludeversion {svntags}

\newtheorem {theorem} {Theorem}

\newtheorem {lemma} [theorem] {Lemma}
\newtheorem {proposition} [theorem] {Proposition}
\newtheorem {corollary} [theorem] {Corollary}

\newcommand {\apclass} [1] {\ensuremath{\mathrm A_{#1}}}

\newcommand {\lclass} [2] {\ensuremath{\mathrm L_{#1} \left( #2 \right) }}
\newcommand {\lsclass} [1] {\ensuremath{\mathit l^{#1} }}

\newcommand {\lclassg} [1] {\ensuremath{\mathrm L_{#1}}}

\newcommand {\hclassg} [1] {\ensuremath{\mathrm H_{#1}}}

\newcommand {\BMO} {\ensuremath {\mathrm {BMO}}}
\newcommand {\VMO} {\ensuremath {\mathrm {VMO}}}

\DeclareMathOperator* {\essinf} {ess\,inf}
\DeclareMathOperator* {\esssup} {ess\,sup}

\excludeversion {comment}

\newcommand {\weightu} {\ensuremath {\mathit u}}

\newcommand {\weightw} {\ensuremath {\mathit w}}

\begin {document}

\title [Complex interpolation of $\apclass {1}$-regular lattices]
{Complex interpolation
\\
of couple $(\mathbf X, \mathbf{BMO})$
\\
for $\mathbf A_1$-regular lattices}
\author {D.~V.~Rutsky}
\email {rutsky@pdmi.ras.ru}
\date {\today}
\address {St.Petersburg Department
of Steklov Mathematical Institute RAS
27, Fontanka
191023 St.Petersburg
Russia}


\begin {abstract}
Recent results of A.~Lerner concerning certain properties of the Fefferman-Stein maximal function
are applied to show 
that $(\BMO, X)_\theta = X^\theta$, $0 < \theta < 1$, for a Banach lattice $X$ of measurable functions on $\mathbb R^n$
satisfying the Fatou property
such that 
$X$ has order continuous norm and the Hardy-Littlewood maximal operator $M$ is bounded in $(X^\alpha)'$
for some $0 < \alpha \leqslant 1$.
\end {abstract}

\maketitle

\setcounter {section} {-1}

\section {Introduction}

\label {anintroduction}
Recently various classical results of harmonic analysis
for important classical Banach spaces such as $\lclassg {p}$ have been generalized to their variable exponent analogues
such as $\lclassg {p (\cdot)}$ and in some cases to general Banach lattices.
Interpolation of such spaces has also received some attention; see, e.~g., \cite {varpbook},
\cite {kempkavybiral2012}, \cite {almeidahasto2012}, \cite {kopaliani2009}.
In particular, in \cite {kopaliani2009} it was established with the help of variable exponent Triebel-Lizorkin spaces
that $(\lclassg {p (\cdot)}, \BMO)_\theta = \lclassg {\frac {p (\cdot)} {1 - \theta}}$
on $\mathbb R^n$ for $0 < \theta < 1$
along with the corresponding formula for $\hclassg {1}$
under the assumption that the Hardy-Littlewood maximal operator $M$ is bounded in $\lclassg {p (\cdot)}$.
This extends the classical result going back to
\cite {feffermanstein1972} saying that in the scale of complex interpolation spaces $\lclassg {p}$
one can replace the endpoint space $\lclassg {\infty}$ by $\BMO$.
In this short note we establish an extension of this result to fairly general Banach lattices.
Although it appears feasible to extend the approach of \cite {kopaliani2009} to this generality by studying the Triebel-Lizorkin type spaces
corresponding to general Banach lattices,
in this case it feels more natural to use a straightforward extension of the original argument involving application of the Fefferman-Stein maximal function,
which is made possible by recent results of A.~Lerner \cite {lerner2010} extending certain properties of the Fefferman-Stein maximal function
to fairly general Banach lattices.
There are, of course, a number of technical difficulties to be addressed.

\section {Preliminaries}

\label {preliminaries}

First, let us define the complex interpolation spaces.
For a couple $(X_0, X_1)$ of compatible complex Banach spaces and $0 \leqslant \theta \leqslant 1$ the complex interpolation space
$(X_0, X_1)_\theta$ is defined as follows (for more detail see, e.~g., \cite [Chapter 4] {bergh}).
Let $\mathcal F_{X_0, X_1}$ be the space of all
bounded and continuous functions $f : z \mapsto f_z$ that are
defined on the strip $S = \{ z \in \mathbb C \mid 0 \leqslant \Re z \leqslant 1\}$ and take values in $X_0 + X_1$
such that $f$ are analytic on the interior of $S$, $f_{j + i t} \in X_j$ for $j \in \{0, 1\}$ and all $t \in \mathbb R$, and
$\|f_{j + i t}\|_{X_j} \to 0$ as $|t| \to \infty$.
The space $\mathcal F_{X_0, X_1}$ is equipped with the norm $\|f\|_{\mathcal F_{X_0, X_1}} = \sup_{t \in \mathbb R, j \in \{0, 1\}} \|f_{j + i t}\|_{X_j}$.
Then space $(X_0, X_1)_\theta = \left\{f_\theta \mid f \in \mathcal F_{X_0, X_1}\right\}$ equipped with the norm
$$
\|a\|_{(X_0, X_1)_\theta} = \inf \left\{\|f\|_{\mathcal F_{X_0, X_1}} \mid f \in \mathcal F, f_\theta = a\right\}
$$
is an interpolation space of exponent $\theta$ between $X_0$ and $X_1$.
Moreover, $X_0 \cap X_1$ is dense in $(X_0, X_1)_\theta$ for $0 < \theta < 1$
(see, e.~g., \cite [Theorem~4.2.2] {bergh}),
and if $X_0 \cap X_1$ is dense in $X_j$ for $j \in \{0, 1\}$ then
$(X_0, X_1)_j = X_j$
(see, e.~g., remarks after \cite [Chapter~4, Theorem~1.3] {kps}).

We are now going to list some well-known standard facts about Banach lattices of measurable functions that we need in the present work;
for more detail see, e.~g., \cite {kantorovichold}.
A Banach space $X$ of measurable functions on a $\sigma$-finite
measurable space $\Omega$ (for example, $\Omega = \mathbb R^n$ with the
Lebesgue measure) is called a \emph {Banach lattice} if for any $f \in X$ and a measurable function $g$ such that $|g| \leqslant f$ almost everywhere
we also have $g \in X$ and $\|g\|_X \leqslant C \|f\|_X$ with some $C$ independent of $f$ and $g$.
We say that $X$ satisfies the \emph {Fatou property}
(which is usually assumed in the literature, implicitly or otherwise)
if $f_n \in X$, $\|f_n\|_X \leqslant 1$ and $f_n \to f$ almost everywhere for some $f$ imply that $f \in X$ and $\|f\|_X \leqslant 1$.
The \emph {order dual} $X'$ of $X$ can be identified with the Banach lattice of measurable functions $g$ having finite norm
$\|g\|_{X'} = \sup_{f \in X, \|f\|_X \leqslant 1} \int_{\Omega} f g$.
The Fatou property of a lattice $X$ is equivalent to order reflexivity of $X$, that is to the relation $X = X''$.
A Banach lattice is said to have an \emph {order continuous norm} if $\|f_n\|_X \to 0$ for every nonincreasing sequence of functions $f_n \in X$ convegring
to $0$ almost everywhere.  A Banach space has order continuous norm if and only if its order dual is isomorphic to the dual Banach space, i.~e.
$X^* = X'$.  Thus, for example, $\lclassg {p}' = \lclassg {p'}$ for $1 \leqslant p \leqslant \infty$ with $\frac 1 p + \frac 1 {p'} = 1$, but
$\lclassg {p}^* = \lclassg {p'}$ holds true only for $1 \leqslant p < \infty$.

For Banach lattices $X_0$, $X_1$ and $0 < \theta < 1$ the \emph {Calderon product}
is the lattice of measurable functions $f$ such that the norm
$$
\|f\|_{X_0^{1 - \theta} X_1^{\theta}} =
\inf \left\{ \left\| |f_0|^{\frac 1 {1 - \theta}} \right\|_{X_0}^{1 - \theta} \left\| |f_1|^{\frac 1 \theta} \right\|_{X_1}^\theta \mid
f = f_0 f_1 \right\}
$$
is finite.  It is well known (see, e.~g., \cite {calderon1964}, \cite {lozanovsky1969}) that if $X_0$ and $X_1$ have the Fatou property then
$X_0^{1 - \theta} X_1^{\theta}$ is also a Banach lattice satisfying the Fatou property and
its order dual can be computed as $\left(X_0^{1 - \theta} X_1^\theta \right)' = X_0'^{1 - \theta} X_1'^\theta$.
Setting $X_0 = \lclassg {\infty}$ and $X^\theta = X^\theta \lclassg {\infty}^{1 - \theta}$ allows one to scale lattices,
so that, for example, $\left[\lclassg {p}\right]^\theta = \lclassg {\frac p \theta}$,
and we have a useful duality relation $\left(X^\theta\right)' = {X'}^\theta \lclassg {1}^{1 - \theta}$.
It is easy to see that if either $X_0$ or $X_1$ has order continuous norm then $X_0^{1 - \theta} X_1^{\theta}$ also has order continuous norm.
In~\cite {calderon1964} (see also \cite [Chapter~4, Theorem~1.14] {kps})
it was established that Calderon products describe complex interpolation spaces between Banach lattices, i.~e.
$(X_0, X_1)_\theta = X_0^{1 - \theta} X_1^\theta$, provided that $X_0^{1 - \theta} X_1^\theta$ has order continuous norm\footnote{
In \cite {kopaliani2009} and in some other papers
it was claimed that $(X_0, X_1)_\theta = X_0^{1 - \theta} X_1^\theta$ when $X_0^{1 - \theta} X_1^\theta$ has the Fatou property.
However, in general the Fatou property only gives $(X_0, X_1)^\theta = X_0^{1 - \theta} X_1^\theta$, and a simple example of two weighted
spaces $\lclass {\infty} {\weightw}$ shows that sometimes $(X_0, X_1)^\theta = X_0^{1 - \theta} X_1^\theta {\raisebox {-0.6mm} {$\supsetneqq$}} (X_0, X_1)_\theta$
in this case.
See, e.~g., \cite [\S13.6] {calderon1964}.
}.

Let $X$ be a Banach lattice of measurable functions on $\Omega$.  The lattice $X (\lsclass {\infty})$
is the space of all measurable functions $f = \{f_j\}_{j \in \mathbb Z}$ on $\Omega \times \mathbb Z$ such that the norm
$\|f\|_{X (\lsclass {\infty})} = \left\|\sup_j |f_j| \right\|_X$ is finite.
This is a particular case of the general construction of a lattice with mixed norm that we will use in the present work.
It is easy to see that if $X$ satisfies the Fatou property then so does $X (\lsclass {\infty})$ and
$\left[ X (\lsclass {\infty}) \right]^\theta = X^\theta (\lsclass {\infty})$ for all $0 < \theta < 1$.
Observe that $X (\lsclass {\infty})$ never has order continuous norm.  Because of this
we will need the following simple proposition (which actually holds true for any lattice of measurable functions in place of $\lsclass {\infty}$;
although we will only use the well-known inclusion $\subset$,
we also prove the converse inclusion for completeness).
\begin {proposition}
\label {infint}
Let $X_0$ and $X_1$ be Banach lattices of measurable functions.  If $X_0$ has order continuous norm
then
$$
\left(X_0 (\lsclass {\infty}), X_1 (\lsclass {\infty}) \right)_\theta = X_0^{1 - \theta} X_1^\theta (\lsclass {\infty}).
$$
\end {proposition}
Indeed, inclusion
$\left(X_0 (\lsclass {\infty}), X_1 (\lsclass {\infty}) \right)_\theta \subset
X_0^{1 - \theta} X_1^\theta (\lsclass {\infty})
$
follows at once from \cite [\S13.6, i] {calderon1964}, and we only need to establish the converse inclusion.
Let $f = \{f_j\}_{j \in \mathbb Z} \in X_0^{1 - \theta} X_1^\theta (\lsclass {\infty})$, and $F = \sup_j |f_j|$.
Then $F \in X_0^{1 - \theta} X_1^\theta = \left(X_0, X_1\right)_\theta$.  
This means that $F = f_\theta$ for some
$f \in \mathcal F_{X_0, X_1}$ with an appropriate estimate on the norm.
Defining $F_z = \{f_{z, j}\}_{\j \in \mathbb Z}$
by $f_{z, j} = \frac {f_j} F f_z$ (with the usual convention that $\frac 0 0 = 0$) shows that
$F_z \in \mathcal F_{X_0 (\lsclass {\infty}), X_1 (\lsclass {\infty})}$ with the same norm as $f_z$, so
$$
f = F_\theta \in \left(X_0 (\lsclass {\infty}), X_1 (\lsclass {\infty}) \right)_\theta
$$
with an appropriate estimate on the norm.
The proof of Proposition~\ref {infint} is complete.

The \emph {Hardy-Littlewood maximal operator $M$}
is defined for all locally summable functions $f$ by
$$
M f (x) = \sup_{Q \ni x} \frac 1 {|Q|} \int_Q |f (y)| dy, \quad x \in \mathbb R^n,
$$
where the supremum is taken over all cubes $Q \subset \mathbb R^n$ containing $x$ with edges parallel to the coordinate axes.
A locally summable nonnegative function $\weightw$ belongs to the
\emph {Muckenhoupt class $\apclass {1}$ with constant $c$} if $M \weightw \leqslant c \weightw$ almost everywhere.
We say that a Banach lattice $X$ of measurable functions on $\mathbb R^n$ is \emph {$\apclass {1}$-regular with constants $(c, m)$}
if for any $f \in X$ there exists some majorant $\weightw \geqslant |f|$ belonging to $\apclass {1}$ with constant $c$ such that
$\|\weightw\|_X \leqslant c \|f\|_X$.
By \cite [Proposition~1.2] {rutsky2011en} a Banach lattice $X$ is $\apclass {1}$-regular if and only if $M$ is bounded in $X$;
thus $\apclass {1}$-regularity of $X$ can justifiably be considered a rather convenient term for boundedness of $M$ in $X$.
The proof is very simple: an $\apclass {1}$-majorant for $f \in X$ gives at once the necessary estimate for $M f$, and conversely
an $\apclass {1}$-majorant can be quickly obtained from boundedness of $M$ in $X$ by the well-known construction due to Rubio de Francia.

With the help of the theory of Muckenhoupt weights it is rather easy to see that
the $\apclass {1}$-regularity property is
``almost self-dual'' in the following sense.
\begin {proposition} {{\cite [Proposition~1.7] {rutsky2011en}}}
\label {a1duality}
Let $X$ be a Banach lattice of measurable functions on $\mathbb R^n$ having either the Fatou property or order continuous norm.
Suppose that $X'$ is $\apclass {1}$-regular.  Then $X^\theta$ is also $\apclass {1}$-regular for any $0 < \theta < 1$.
\end {proposition}

The following well-known characterization of $\apclass {1}$ weights is very useful;
in can be found in, e.~g., \cite [Chapter~5, \S5.2] {stein1993}.
\begin {proposition}
\label {maxfunpot}
Let $\weightw$ be a nonnegative locally summable function.
Then $\weightw \in \apclass {1}$ with a constant $c$ if and only if
there exists an exponent $0 < q < 1$, a locally summable function $f$ and constants $c_0, c_1 > 0$
such that $c_0 \weightw \leqslant (M f)^q \leqslant c_1 \weightw$.
If this holds true then constant $c$ and constants $q$, $c_0$, $c_1$ can be estimated in terms of one another.
\end {proposition}
Proposition~\ref {maxfunpot} is a consequence of the reverse H\"older inequality satisfied by $\apclass {1}$ weights.
It allows a very easy proof of the following result.
\begin {proposition}
\label {a1interp}
Let $X$ be an $\apclass {1}$-regular Banach lattice of measurable functions on $\mathbb R^n$.
Then lattices $X^\theta$ and 
$X^{1 - \theta} \lclassg {1}^\theta$ are also $\apclass {1}$-regular for all $0 < \theta < 1$.
\end {proposition}
Indeed, $\apclass {1}$-regularity of $X^\theta$
is a trivial corollary to Proposition~\ref {maxfunpot}, and it is otherwise established at once using the H\"older inequality.
More generally, the H\"older inequality shows that for any two $\apclass {1}$-regular lattices $A$ and $B$ lattice $A^{1- \theta} B^\theta$
is also $\apclass {1}$-regular (see, e.~g., \cite [Proposition~3.4] {rutsky2011en}), and this also implies $\apclass {1}$-regularity of $X^\theta$
since lattice $\lclassg {\infty}$ is trivially $\apclass {1}$-regular.
It is, however, well known that lattice $\lclassg {1}$ is not $\apclass {1}$-regular,
so $\apclass {1}$-regularity of $X^{1 - \theta} \lclassg {1}^\theta$ is a bit more tricky.
Suppose that $f \in X^{1 - \theta} \lclassg {1}^\theta$.
We may assume that $f \geqslant 0$ and $\|f\|_{X^{1 - \theta} \lclassg {1}^\theta} = 1$.
Then $f = g^{1 - \theta} h^\theta$ with some $g \in X$ and $h \in \lclassg {1}$ with $\|g\|_X \leqslant 2$ and $\|h\|_{\lclassg {1}} \leqslant 2$.
Let $\weightw$ be an $\apclass {1}$-majorant for $g$ in $X$.  Then by Proposition~\ref {maxfunpot} weight
$\weightw$ is pointwise equivalent to $(M a)^q$ almost everywhere with some locally summable function $a$ and with $0 < q < 1$
depending only on the $\apclass {1}$-regularity constants of $X$.
Since $M$ is bounded in $\lclassg {p}$ for any $1 < p \leqslant \infty$ we have an estimate
$\left\|(M [h^\alpha])^{\frac 1 \alpha}\right\|_{\lclassg {1}} \leqslant c$ with some $c$ independent of $f$ for any $0 < \alpha < 1$.
Observe that $f$ is dominated by $\weightu = c_1 (M a)^{q (1 - \theta)} (M h^\alpha)^{\frac 1 \alpha \theta}$
and $\|\weightu\|_{X^{1 - \theta} \lclassg {1}^\theta} \leqslant c_2$ with some $c_1$ and $c_2$ independent of $f$.
We claim that with a certain choice of $\alpha$ we have $\weightu \in \apclass {1}$ with a constant independent of $f$.
Indeed, by the H\"older inequality
\begin {equation}
\label {qconve}
\frac 1 {c_1 |Q|} \int_Q \weightu \leqslant
\left(\frac 1 {|Q|} \int_Q (M a)^{p q (1 - \theta)}\right)^{\frac 1 p}
\left(\frac 1 {|Q|} \int_Q \left(M h^\alpha\right)^{\frac 1 \alpha p' \theta}\right)^{\frac 1 {p'}}
\end {equation}
for any cube $Q \subset \mathbb R^n$ and $1 < p < \infty$.  If we choose the parameters so that
\begin {equation}
\label {pq1t}
p q (1 - \theta) < 1, \quad
\frac {p'} \alpha \theta < 1,
\end {equation}
then $(M a)^{p q (1 - \theta)} \in \apclass {1}$ and $\left(M h^\alpha\right)^{\frac 1 \alpha p' \theta} \in \apclass {1}$
by Proposition~\ref {maxfunpot}, and therefore \eqref {qconve} and \eqref {pq1t} imply that
$$
\frac 1 {c_1 |Q|} \int_Q \weightu \leqslant c_2 (M a)^{q (1 - \theta)} (M h^\alpha)^{\frac 1 \alpha \theta} = c_2 \weightu
$$
almost everywhere with some constant $c_2$,
i.~e. $\weightu \in \apclass {1}$ with an appropriate estimate on the constant.
Rewriting \eqref {pq1t}
as $\frac \alpha {\alpha - \theta} < p < \frac 1 {q (1 - \theta)}$
we see that we can always choose an appropriate $p$
if we take any $1 > \alpha > \frac \theta {1 - q (1 - \theta)}$.  The proof of Proposition~\ref {a1interp} is complete.

Let $f$ be a measurable function on $\mathbb R^n$. The \emph { nonincreasing rearrangement $f^*$} of $f$ is defined by
$$
f^* (t) = \inf \left\{ \lambda > 0 \mid |\{ x \in \mathbb R^n \mid |f (x)| > \lambda \}| \leqslant t \right\}, \quad 0 < t < \infty.
$$
Let $S_0$ be the set of all measurable functions $f$
on $\mathbb R^n$ such that
$$
f^* (+\infty) = \lim_{t \to \infty} f^* (t) = 0.
$$
It is easy to see that $S_0$ contains all measurable functions supported on sets of finite measure
and also $\lclassg {p} \subset S_0$ for all $0 < p < \infty$.
Thus if $X$ is a Banach lattice of measurable functions having order continuous norm then $S_0 \cap X$ is dense in $X$.
Density of $S_0 \cap X$ in a lattice $X$ is a somewhat more general assumption than density of simple functions with compact support in $X$;
for example, simple functions with compact support are not dense in
a lattice $\lclass {\infty} {\weightw} = \weightw \lclassg {\infty}$ with weight $\weightw (x) = (1 + |x|)^{-1}$ but at the same time
we have $\lclass {\infty} {\weightw} \subset S_0$.
\begin {comment}
\begin {proposition}
\label {acns0}
Suppose that $X$ is a Banach lattice of measurable functions on a $\sigma$-finite measurable space $\Omega$
having an order continuous norm.  Then $S_0 \cap X$ is dense in $X$.
\end {proposition}
Indeed, let $f \in X$ and let
$A_n \subset \Omega$ be an increasing sequence of sets of finite measure such that $\bigcup_n A_n = \Omega$.
We set $f_n = f \chi_{A_n}$.  Then $f_n \in S_0$, because $f_n$ are supported on sets $A_n$ of finite measure, 
and $\|f - f_n\|_X \to 0$ by the order continuity of the norm of $X$ since $f - f_n = f_n \chi_{\Omega \setminus A_n}$
is monotonely decreasing and converges to $0$ everywhere.
\end {comment}



Now we will briefly discuss some of the results involving the Fefferman-Stein sharp maximal
function. 
The Fefferman-Stein maximal function $f^\sharp$ on $\mathbb R^n$ is defined for a locally integrable function $f$ by
$$
f^\sharp (x) = \sup_{Q \ni x} \frac 1 {|Q|} \int_Q \left|f (y) - f_Q \right| dy, \quad x \in \mathbb R^n,
$$
where the supremum is taken over all cubes $Q \subset \mathbb R^n$ 
containing $x$ 
with edges parallel to the coordinate axes and and $f_Q = \frac 1 {|Q|} \int_Q f (z) dz$ is the average value of $f$
over $Q$ with respect to the Lebesgue measure.
Space $\BMO$ can then be defined as the space of all locally integrable functions $f$ such that $f^\sharp \in \lclassg {\infty}$ modulo constants
equipped with the norm $\|f\|_{\BMO} = \|f^\sharp\|_{\lclassg {\infty}}$ that turns $\BMO$ into a Banach space; for more detail see, e.~g.,
\cite [Chapter 4] {stein1993}.  We have continuous inclusion $\lclassg {\infty} \subset \BMO$, but $\lclassg {\infty}$ is a proper 
subspace of $\BMO$.  The usefulness of $\BMO$ in harmonic analysis stems mainly from the fact that $\BMO$ is dual to the Hardy space $\hclassg {1}$ and
many operators of interest are not bounded in $\lclassg {\infty}$ but act boundedly from
$\lclassg {\infty}$ to $\BMO$ if suitably defined on this space.

\begin {comment}
It was observed in \cite {jawerthtorchinsky1985} that the real-valued $\BMO$ can be characterized by the \emph {Str\"omberg local sharp maximal function}
\begin {equation}
\label {sharpfdef}
M^\sharp_\lambda f (x) = \sup_{Q \ni x} \inf_{c \in \mathbb R} ((f - c)\chi_Q)^* (\lambda |Q|), \quad x \in \mathbb R^n,
\end {equation}
where $f \in S_0$ and the supremum is taken over all cubes $Q \subset \mathbb R^n$ containing the point $x$.
This is due to the following estimate (see, e.~g., \cite {jawerthtorchinsky1985}, \cite {lerner2003})
that holds true
with all sufficiently small $0 < \lambda < 1$ and some $c_0, c_1 > 1$
for all locally summable $f$:
\begin {equation}
\label {sharpequiv}
c_0 M M^\sharp_\lambda f (x) \leqslant f^\sharp (x) \leqslant c_1 M M^\sharp_\lambda f (x), \quad x \in \mathbb R^n.
\end {equation}
\end {comment}
\begin {theorem} [{\cite [Corollary~4.3] {lerner2010}}]
\label {fsharpdual}
Suppose that $X$ is an $\apclass {1}$-regular real Banach lattice of measurable functions on $\mathbb R^n$
having the Fatou property.  Then the following conditions are equivalent.
\begin {enumerate}
\item
$X'$ is $\apclass {1}$-regular.
\item
There exists some $c > 0$ such that $\|f\|_X \leqslant c \|f^\sharp\|_X$ for all $f \in S_0 \cap X$.
\end {enumerate}
\end {theorem}
This theorem can be considered an extension of well-known classical results for $X = \lclassg {p}$ (see, e.~g., \cite {stein1993}).
The proof involves a certain linearization of $M$, 
pointwise equivalence of $f^\sharp$ and $M M_\lambda^\sharp$ for some $\lambda$ (where $M^\sharp_\lambda$ denotes the Str\"omberg local
sharp maximal function) and
the fact that $M_\lambda^\sharp$ is dual to $M$ in the sense that $\int |f g| \leqslant c \int M_\lambda^\sharp f\, M g$
for suitable $f$ and $g$.

It is easy to see
that the estimate in Theorem~\ref {fsharpdual} can be extended to the entire lattice $X$ provided that $S_0 \cap X$ is dense in $X$,
and the complex lattices can be included as well.
\begin {proposition}
\label {fsharpest}
Suppose that $X$ is a Banach lattice of measurable functions on $\mathbb R^n$ having the Fatou property,
both $X$ and $X'$ are $\apclass {1}$-regular
and $S_0 \cap X$ is dense in $X$.  Then there exists some $c > 0$ such that
\begin {equation}
\label {eqfsh}
\|f\|_X \leqslant c \|f^\sharp\|_X
\end {equation}
for all $f \in X$.
\end {proposition}

Indeed, suppose that $f \in X$ under the conditions of Proposition~\ref {fsharpest}, $f$ is real,
and let $f_n \in S_0 \cap X$ be a sequence such that $f_n \to f$ in $X$.
Observe that the Fefferman-Stein maximal function is subadditive and
$g^\sharp \leqslant 2 M g$ for all locally summable functions $g$.  Therefore Theorem~\ref {fsharpdual} allows us to carry out the estimate
\begin {multline}
\label {eqfsha}
\frac 1 c \|f_n\|_X \leqslant \|f_n^\sharp\|_X \leqslant \|f^\sharp\|_X + \|(f - f_n)^\sharp\|_X \leqslant
\\
\|f^\sharp\|_X + 2 \|M (f - f_n)\|_X \leqslant \|f^\sharp\|_X + 2 \|M\|_{X \to X} \|f - f_n\|_X.
\end {multline}
Passing to the limit $n \to \infty$ in \eqref {eqfsha} yields \eqref {eqfsh} for all real functions $f \in X$.
If $f \in X$ is complex then \eqref {eqfsh} implies that
$\|\Re f\|_X \leqslant c \|(\Re f)^\sharp\|_X \leqslant c \|f^\sharp\|_X$ because
$(\Re f)^\sharp \leqslant f^\sharp$ almost everywhere, and the same estimate holds true for $\Im f$.
Combining these estimates together yields
$$
\|f\|_X \leqslant \|\Re f\|_X + \|\Im f\|_X \leqslant 2 c \|f^\sharp\|_X.
$$

\begin {comment}
It was observed in \cite {lerner2010} that the estimate \eqref {eqfsh} is self-improving in a certain indirect way.
It is not hard to see with the help of \eqref {sharpequiv} and a few simple properties of the local sharp maximal function
that
\eqref {eqfsh} is self-improving in a rather direct way as well.
A measurable function $u$ is called \emph {unimodular} if $|u| = 1$ almost everywhere.
\begin {proposition}
\label {naxmodulus}
Let $u$ be a unimodular function, $0 < \lambda < 1$ and $f \in S_0$.
Then
$M^\sharp_\lambda f \leqslant M^\sharp_{\frac \lambda 2} (f u)$ almost everywhere.
\end {proposition}
Indeed, let $Q \subset \mathbb R^n$ be a cube and take some $\varepsilon > 0$.
There exists some $c \in \mathbb R$ such that
\begin {equation}
\label {ccontr}
(|f u - c| \chi_Q)^* \left(\frac \lambda 2 |Q|\right) < \inf_{c_1 \in \mathbb R} (|f u - c_1| \chi_Q)^* \left(\frac \lambda 2 |Q|\right) + \varepsilon.
\end {equation}
Let
$$
E = \left\{ x \in Q \mid |f u - c| > (|f u - c| \chi_Q)^* \left(\frac \lambda 2 |Q|\right) + \varepsilon \right\}.
$$
Then $|E| \leqslant \frac \lambda 2 |Q|$ by the definition of the nonincreasing rearrangement.
Let $E_\sigma = \{x \in E \mid u (x) = \sigma\}$ for $\sigma \in \{-1, 1\}$.  Then $E$ is
the disjoint union of $E_\sigma$ up to a set of measure $0$.
Therefore for some $\sigma$.

\begin {theorem}
\label {sharpsi}
Let $X$ be a
real Banach lattice of measurable functions on $\mathbb R^n$, and suppose that
\begin {equation}
\label {fsh3}
\|f\|_X \leqslant c \|f^\sharp\|_X
\end {equation}
for all $f \in X$ with some $c > 0$ independent of $f$.
Then for all sufficiently small $p > 1$ and unimodular functions $u$ we also have
\begin {equation}
\label {eqrshi}
\|f\|_X \leqslant c_p \left\|\left(\left[|f|^{p} u\right]^\sharp\right)^{\frac 1 p}\right\|_X
\end {equation}
for all $f \in X$ with some $c_p > 0$ independent of $f$ and $u$, and also
\begin {equation}
\label {eqrshith}
\|f\|_{X^\theta} \leqslant c_\theta \|f^\sharp\|_{X^\theta}
\end {equation}
for all $f \in X^\theta$ for all $\theta_0 < \theta \leqslant 1$ with some $0 < \theta_0 < 1$ and $c_\theta > 0$.
\end {theorem}
Of course, \eqref {eqrshi} implies \eqref {eqrshith}.
Observe that the $\apclass {1}$-regularity condition is self
improving in the sense that there exists some $r > 1$ such that $X$ is locally $\lclassg {r}$-summable and the scaled maximal operator
$$
M_r f (x) = \left(\sup_{Q \ni x} \frac 1 {|Q|} \int_Q |f (y)|^r dy\right)^{\frac 1 r}, \quad x \in \mathbb R^n,
$$
is bounded in $X$.  See \cite {lernerperez2007} for a rather general approach; it should be noted, however,
that it is an easy consequence of the characterization of $\apclass {1}$-regularity in terms of $\apclass {1}$-weights
(see, e.~g., \cite [Proposition~1.2] {rutsky2011en}) and their self-improving property (see, e.~g., \cite [Chapter~5, \S5.2] {stein1993}).
The H\"older inequality implies that $M f \leqslant M_r f$ almost everywhere for suitable $f$.
Making use of \eqref {sharpequiv}, we get from \eqref {fsh3} an estimate
\end {comment}

It is easy to see that $\lclassg {\infty}^\alpha \subset \lclassg {\infty}$ for all $0 < \alpha \leqslant 1$.
It is also not hard to verify that $\BMO$ (which is also a lattice) satisfies the same property.
\begin {proposition}
\label {bmofold}
Suppose that $f \in \BMO$ and $f \geqslant 0$.
Then $f^\alpha \in \BMO$ for all $0 < \alpha \leqslant 1$.
\end {proposition}
Since  $f^\alpha - f^\alpha \vee 1$ is a bounded function under the conditions of
Proposition~\ref {bmofold} (and hence $f^\alpha - f^\alpha \vee 1 \in \BMO$),
it suffices to verify that $f^\alpha \vee 1 = (f \vee 1)^\alpha \in \BMO$.
We have $g = f \vee 1 \in \BMO$ because $\BMO$ is a lattice,
and then $g^\alpha \in \BMO$ is clear because the map $F : y \mapsto y^\alpha$ is contractive for
$y \geqslant 1$ and therefore oscillations of $g^\alpha$ do not increase compared to the corresponding oscillations of $g$.
Perhaps the easiest way to verify this formally is via the Str\"omberg characterization of $\BMO$ mentioned above
(see, e.~g., \cite [Chapter~4, \S6.6] {stein1993}; this also involves the local sharp maximal function $M_\lambda^\sharp$
used in the proof of Theorem~\ref {fsharpdual})
which states that $g \in \BMO$ if and only if
there exist some constants $0 < \gamma < \frac 1 2$ and $\lambda > 0$ such that
\begin {equation}
\label {gstrom}
\inf_{c_Q \in \mathbb R} |\{x \in Q \mid |g (x) - c_Q| > \lambda \}| \leqslant \gamma |Q|
\end {equation}
for any cube $Q \subset \mathbb R^n$.  It is easy to see that if $F$ is a contractive map then \eqref {gstrom}
implies 
\begin {equation*}
\inf_{c_Q \in \mathbb R} |\{x \in Q \mid |F \circ g (x) - F (c_Q)| > \lambda \}| \leqslant \gamma |Q|
\end {equation*}
for any cube $Q \subset \mathbb R^n$, so $F \circ g \in \BMO$ if $g \in \BMO$.

\begin {comment}
Another way to see this, which also shows that the restriction
$\alpha \leqslant 1$ is superfluous,
is to observe that $y \mapsto 1 \vee y^\alpha$ is a Lipschitz map and
$\BMO$ is stable under composition on the left with Lipschitz maps,
which is easy to verify using the Str\"omberg characterization of $\BMO$ mentioned above
(see, e.~g., \cite [Chapter~4, \S6.6] {stein1993}).
\end {comment}

\begin {proposition}
\label {bmoadropi}
Let $X$ be a Banach lattice and suppose that
$X^\alpha \cap \BMO$ is a subset of $X^{\theta \alpha}$ for some $0 < \alpha, \theta < 1$.
Then $X \cap \BMO$ is a subspace of $X^\theta$ for all $0 < \eta < 1$.
\end {proposition}
Indeed, since $\BMO$ is a lattice, it is sufficient verify the inclusion $X \cap \BMO \subset X^\eta$ for
nonnegative functions.  Suppose that $f \in X \cap \BMO$ and $f \geqslant 0$ almost everywhere.
Then $f^\alpha \in X^\alpha$ and by Proposition~\ref {bmofold} we have $f^\alpha \in \BMO$.
Thus $f^\alpha \in X^\alpha \cap \BMO \subset X^{\theta \alpha}$
and therefore $f \in X^\theta$.

\begin {comment}
Perhaps the easiest way to do this is via the Str\"omberg characterization of $\BMO$ mentioned above
(see, e.~g., \cite [Chapter~4, \S6.6] {stein1993}) which states that $g \in \BMO$ if and only if
there exist some constants $0 < \gamma < \frac 1 2$ and $\lambda > 0$ such that
\begin {equation}
\label {gstrom}
\inf_{c_Q \in \mathbb R} |\{x \in Q \mid |g (x) - c_Q| > \lambda \}| \leqslant \gamma |Q|
\end {equation}
for any cube $Q \subset \mathbb R^n$.
Indeed, if 
\end {comment}

\begin {comment}
\begin {proposition}
\label {s0densth}
Let $X$ be a Banach lattice of measurable functions such that $S_0 \cap X$ is dense in $X$.
Then $S_0 \cap X^\theta$ is dense in $X^\theta$ for all $0 < \theta < 1$.
\end {proposition}
It is sufficient to establish that $S_0 \cap X^\theta_+$ is dense in
$$
X^\theta_+ = \{f \in X^\theta \mid f \geqslant 0 \text { a. e.}\}.
$$
Suppose that $f \in X^\theta_+$.  Then $f^{\frac 1 \theta} \in X$, and so
for any $\varepsilon > 0$ there exists
some $g \in S_0 \cap X_+$ such that $\left\| f^{\frac 1 \theta} - g \right\|_X \leqslant \varepsilon$;
by replacing $g$ by $g \wedge f^{\frac 1 \theta}$
we may additionally assume that $g \leqslant f^{\frac 1 \theta}$.
By the well known inequality $(a + b)^\theta \leqslant a^\theta + b^\theta$ for nonnegative $a$ and $b$ we have
$$
f = \left(f^{\frac 1 \theta} \right)^\theta \leqslant \left( f^{\frac 1 \theta} - g \right)^\theta + g^\theta,
$$
so
$$
\left\|f - g^\theta\right\|_{X^\theta}^{\frac 1 \theta} \leqslant \left\| f^{\frac 1 \theta} - g \right\|_X \leqslant \varepsilon,
$$
which concludes the proof of Proposition \ref {s0densth}.
\end {comment}

\section {Interpolation}

\label {sinterp}

We are now ready to state the main result.
\begin {theorem}
\label {bintt}
Suppose that $X$ is a Banach lattice of measurable functions on $\mathbb R^n$ having the Fatou property and order continuous norm,
and lattice 
$(X^\alpha)'$ is $\apclass {1}$-regular for some $0 < \alpha \leqslant 1$.
Then
\begin {equation}
\label {bmointerp}
(\BMO, X)_\theta = X^\theta
\end {equation}
for any $0 < \theta < 1$.
\end {theorem}

We will give a few remarks before passing to the proof of Theorem~\ref {bintt}.
The assumption that $(X^\alpha)'$ is $\apclass {1}$-regular
combined with the assumption that $X$ has order continuous norm cannot be dropped from Theorem~\ref {bintt}.
Otherwise we would have had
\begin {equation}
\label {bmobmogus}
(\BMO_b, \lclassg {\infty})_\theta = (\BMO, \lclassg {\infty})_\theta = \lclassg {\infty},
\end {equation}
where $\BMO_b$ is the closure of $\lclassg {\infty}$ in $\BMO$.
Equation \eqref {bmobmogus} implies that $\BMO_b = \lclassg {\infty}$ by \cite [Theorem~1.7] {stanfey1970}.
However, it is well known that
$$
\lclassg {\infty} \neq \VMO \subset \BMO_b;
$$
see, e.~g., \cite [Chapter~4, \S6.8] {stein1993}.
On the other hand, it seems that $\apclass {1}$-regularity of $(X^\alpha)'$ should imply order continuity of the norm of $X$.
This is true at least in the case of variable exponent Lebesgue spaces $X = \lclassg {p (\cdot)}$ since
$(X^\alpha)' = \lclassg {\left({p (\cdot)} \slash \alpha\right)'}$ and
$\esssup p (\cdot) = \infty$ would imply $\essinf \left( {p (\cdot)} \slash \alpha\right)' = 1$ which contradicts
$\apclass {1}$-regularity of $(X^\alpha)'$ by \cite [Theorem~4.7.1] {varpbook}.

It is easy to see that if the conditions of Theorem~\ref {bintt} are satisfied for some $\alpha$ then they are satisfied
for all smaller values of $\alpha$, and the lattice $X^\beta$ is $\apclass {1}$-regular for all $0 < \beta < \alpha$.
Indeed, under the conditions of Theorem~\ref {bintt} lattice
$X^\beta = (X^\alpha)^{\frac \beta \alpha}$
is $\apclass {1}$-regular for any
$0 < \beta < \alpha$ by Proposition~\ref {a1duality}, and lattice
$$
(X^\beta)' = (X')^\beta \lclassg {1}^{1 - \beta} = \left[(X')^\alpha \lclassg {1}^{1 - \alpha}\right]^{\frac \beta \alpha} \lclassg {1}^{1 - \frac \beta \alpha} =
[(X^\alpha)']^{\frac \beta \alpha} \lclassg {1}^{1 - \frac \beta \alpha}
$$
is $\apclass {1}$-regular for the same
values of $\beta$ by Proposition~\ref {a1interp}.

We now provide a couple of applications for Theorem~\ref {bintt}.
Muckenhoupt weights $\weightw \in \apclass {p}$ for $1 < p < \infty$ are exactly those
for which the Hardy-Littlewood maximal operator $M$ is bounded in the weighted Lebesgue space $\lclass {p} {\weightw}$ with norm defined by
$$
\|f\|_{\lclass {p} {\weightw}}^p = \int |f|^p \weightw
$$
(here we use this classical definition for the sake of simplicity;
in \cite {rutsky2011en}, for example, the same space was denoted by $\lclass {p} {\weightw^{-\frac 1 p}}$ which gives more consistency with the endpoint
$p = \infty$ and Calderon products).  We can naturally extend this definition to $p = \infty$ by
$\apclass {\infty} = \bigcup_{p > 1} \apclass {p}$; for more detail on Muckenhoupt weights see,
e.~g., \cite [Chapter~5] {stein1993}.
\begin {corollary}
\label {wl1c}
Suppose that $\weightw \in \apclass {\infty}$.  Then for any $0 < \theta < 1$ we have
\begin {equation}
\label {bmoweight}
(\BMO, \lclass {1} {\weightw})_\theta = \lclass {\frac 1 \theta} {\weightw}.
\end {equation}
\end {corollary}
We want to verify that the conditions of Theorem~\ref {bintt} are satisfied for $X = \lclass {1} {\weightw}$ under the conditions of Corollary~\ref {wl1c}.
Indeed,
$$
(\lclass {1} {\weightw}^\alpha)' = (\lclass {\frac 1 \alpha} {\weightw})' = \lclass {\frac 1 {1 - \alpha}} {\weightw^{-1}},
$$
and $\apclass {1}$-regularity of this lattice for suitable values of $\alpha$ follows from the following simple proposition.
\begin {lemma}
\label {ainfinv}
Suppose that $\weightw \in \apclass {\infty}$.  Then $\weightw^{-1} \in \apclass {\infty}$.
\end {lemma}
There are many straightforward
ways to establish Lemma~\ref {ainfinv} using numerous characterizations of the $\apclass {\infty}$ weights; here we are going to use
nothing more than Proposition~\ref {a1duality}.
Indeed, by the assumptions we have $\weightw \in \apclass {p_0}$ with some $1 < p_0 < \infty$. 
Then lattice $\left[\lclass {p'} {\weightw^{-1}}\right]' = \lclass {p} {\weightw}$
is $\apclass {1}$-regular for all $p \geqslant p_0$, and by Proposition~\ref {a1duality}
lattice
$$
\left[\lclass {p'} {\weightw^{-1}} \right]^{\frac 1 q} = \lclass {p' q} {\weightw^{-1}}
$$
is $\apclass {1}$-regular for all $p \geqslant p_0$ and $q > 1$, so $\weightw^{-1} \in \apclass {p' q} \subset \apclass {\infty}$
as claimed.

Application of Theorem~\ref {bintt} to the case $X = \lclassg {p (\cdot)}$ yields part of the results from
\cite {kopaliani2009}; for definitions and general discussion of variable exponent Lebesgue
spaces $\lclassg {p (\cdot)}$ see, e.~g., \cite {varpbook}.
\begin {corollary}
\label {wl1v}
Let $p (\cdot) : \mathbb R^n \to [1, \infty]$ be a measurable function such that $\esssup_{x \in \mathbb R^n} p (x) < \infty$
and suppose that $\lclassg {p (\cdot)}$ is $\apclass {1}$-regular.  Then
\begin {equation}
\label {bmovarp}
(\BMO, \lclassg {p (\cdot)})_\theta = \lclassg {\frac {p (\cdot)} \theta}
\end {equation}
for all $0 < \theta < 1$.
\end {corollary}
Variable exponent Lebesgue spaces $\lclassg {p (\cdot)}$ can be regarded as a natural generalization of the standard Lebesgue spaces
$\lclassg {p}$, which correspond to the case $p (\cdot) = p$,
and lattice operations in spaces $\lclassg {p (\cdot)}$ behave largely the same as their Lebesgue space counterparts.
Observe that $\apclass {1}$-regularity of $\lclassg {p (\cdot)}$ implies by \cite [Theorem~4.7.1] {varpbook} that
$\essinf_{x \in \mathbb R^n} p (x) > 1$, and by \cite [Theorem~5.7.2] {varpbook} it follows that
lattice
$\left[\lclassg {p (\cdot)}\right]' = \lclassg {p' (\cdot)}$ is also $\apclass {1}$-regular.
Condition $\esssup_{x \in \mathbb R^n} p (x) < \infty$ easily implies (see, e.~g., \cite [Lemma~2.3.16] {varpbook})
that $\lclassg {p (\cdot)}$ has order continuous norm.  Therefore we can apply Theorem~\ref {bintt} with $\alpha = 1$,
which concludes the proof of Corollary~\ref {wl1v}.

We now begin the proof of Theorem~\ref {bintt}.  First, observe that T.~Wolff's well-known result concerning gluing of interpolation scales
allows us to reduce it to the case $\alpha = 1$.
\begin {theorem} [{\cite [Theorem~2] {wolff1982}}]
\label {glusc}
Let $A_1$, $A_2$, $A_3$, $A_4$ be Banach spaces.
Suppose that $A_1 \cap A_4$ is a dense subspace of $A_2$ and of $A_3$ and
$$
A_3 = (A_2, A_4)_\gamma, \quad A_2 = (A_1, A_3)_\delta
$$
with some $0 < \gamma, \delta < 1$.
Then
$$
A_2 = (A_1, A_4)_\xi, \quad A_3 = (A_1, A_4)_\psi
$$
for
$\xi = \frac {\gamma \delta} {1 - \delta + \gamma \delta}$ and $\psi = \frac {\gamma} {1 - \delta + \gamma \delta}$.
\end {theorem}
Indeed, suppose that under the conditions of Theorem~\ref {bintt} we have established that
\begin {equation}
\label {alphadrop}
(\BMO, X^\alpha)_\eta = X^{\eta \alpha}
\end {equation}
for all $0 < \eta < 1$.
First, suppose that $\theta < \alpha$ and let $A_1 = \BMO$, $A_2 = X^\theta$, $A_3 = X^\alpha$ and $A_4 = X$.
Equation~\eqref {alphadrop} implies that $\BMO \cap X^\alpha$ is a subspace of $X^{\eta \alpha}$,
so $A_1 \cap A_4$ is a subspace of $A_2$ and $A_3$ by Proposition~\ref {bmoadropi}.
The density assumptions of Theorem~\ref {glusc} are satisfied because
$\BMO \cap X \supset \lclassg {\infty} \cap X$, which is a dense subspace of $(\lclassg {\infty}, X)_\zeta = X^\zeta$ for all
$0 < \zeta < 1$.
The conditions of Theorem~\ref {glusc}
are satisfied with values $\delta = \frac \theta \alpha$ and $\gamma = \frac {\alpha - \theta} {1 - \theta}$,
and thus $X^\theta = A_2 = (A_1, A_4)_\xi = (\BMO, X)_\theta$ ($\xi = \theta$ follows from an easy computation),
i.~e. \eqref {bmointerp} is satisfied for all $0 < \theta < \alpha$; we also get
$X^\alpha = A_3 = (A_1, A_4)_\psi = (\BMO, X)_\alpha$, which is \eqref {bmointerp} for $\theta = \alpha$.
The remaining case $\alpha < \theta < 1$ is then easily established by the reiteration theorem (see, e.~g., \cite [Theorem~4.6.1] {bergh}):
we have
\begin {multline*}
X^\theta = (X^\alpha, X)_\eta = \left( (\BMO, X)_\alpha, (\BMO, X)_1 \right)_\eta =
\\
(\BMO, X)_{(1 - \eta) \alpha + \eta} = (\BMO, X)_\theta
\end {multline*}
for $\eta = \frac {\theta - \alpha} {1 - \alpha}$.

Thus we only need to verify 
\eqref {alphadrop} for all sufficiently small $\alpha$ under the conditions of Theorem~\ref {bintt}.
Since we can always make $\alpha$ smaller, we may assume that lattices $X^\beta$ and $(X^\beta)'$ are $\apclass {1}$-regular for all
$0 < \beta \leqslant \alpha$.
For convenience we replace $X^\alpha$ by $X$; thus lattices $X^\beta$ and $(X^\beta)'$ are $\apclass {1}$-regular for $0 < \beta \leqslant 1$,
and we need to verify that $(\BMO, X)_\eta = X^{\eta}$
for all $0 < \eta < 1$.
The proof now follows the standard pattern.
Let $0 < \theta < 1$.
Since $\lclassg {\infty} \subset \BMO$,
we have $(\BMO, X)_\theta \supset (\lclassg {\infty}, X)_\theta = X^\theta$, and only the converse inclusion needs to be established.
Because $\BMO \cap X$ is dense in $(\BMO, X)_\theta$, it suffices to verify this inclusion on $\BMO \cap X$.
Suppose that $a \in (\BMO, X)_\theta \cap (\BMO \cap X)$; then $a = f_\theta$ with some $f \in \mathcal F_{\BMO, X}$ with
$\|f\|_{\mathcal F_{\BMO, X}} \leqslant 2 \|a\|_{(\BMO, X)_\theta}$.
We enumerate all cubes $\{Q_j\}_{j \in \mathbb N}$ containing $0$ and having rational coordinates of the vertices
and define a function $g = \{g_j\}_{j \in \mathbb N}$ on the strip $S$ by
$$
g_{z, j} (x) = \frac 1 {|Q_j|} \int_{Q_j + x} \left( f_z - \frac 1 {|Q_j|} \int_{Q_j + x} f_z \right)
\frac {\overline {f_\theta - \frac 1 {|Q_j|} \int_{Q_j + x} f_\theta}} {\left|f_\theta - \frac 1 {|Q_j|} \int_{Q_j + x} f_\theta\right|}
$$
for all $z \in S$, $x \in \mathbb R^n$ and $j \in \mathbb N$.
It is easy to see that $g$ is continuous on the strip $S$ and analytic in the interior of $S$.
Moreover, we have estimates
\begin {multline*}
\sup_j |g_{i t, j} (x)| \leqslant 
\sup_j \frac 1 {|Q_j|} \int_{Q_j + x} \left| f_{i t} (x) - \frac 1 {|Q_j|} \int_{Q_j + x} f_{i t} (x) \right| \leqslant
\\
f_{i t}^\sharp (x) \leqslant
\|f_{i t}\|_{\BMO} \leqslant \|f\|_{\mathcal F_{\BMO, X}}
\end {multline*}
and
$$
\sup_j |g_{1 + i t, j} (x)| \leqslant 2 M f_{1 + i t} (x)
$$
for all $t \in \mathbb R$ and almost all $x \in \mathbb R^n$,
so $\|g_{i t}\|_{\lclassg {\infty} (\lsclass {\infty})} \leqslant \|f\|_{\mathcal F_{\BMO, X}}$ and
$\|g_{1 + i t}\|_{X (\lsclass {\infty})} \leqslant 2 \|M f_{1 + i t}\|_X \leqslant c_1 \|f_{1 + i t}\|_X \leqslant c_1 \|f\|_{\mathcal F_{\BMO, X}}$
for all $t \in \mathbb R$ with some constant $c_1 > 1$ independent of $a$.  These estimates also imply that
$\|g_{i t}\|_{\lclassg {\infty} (\lsclass {\infty})} \to 0$ and $\|g_{1 + i t}\|_{X (\lsclass {\infty})} \to 0$ as $t \to \infty$.
Thus
$g \in \mathcal F_{\lclassg {\infty} (\lsclass {\infty}), X (\lsclass {\infty})}$ and
$\|g\|_{\mathcal F_{\lclassg {\infty} (\lsclass {\infty}), X (\lsclass {\infty})}} \leqslant c_1 \|f\|_{\mathcal F_{\BMO, X}} \leqslant 2 c_1 \|a\|_{(\BMO, X)_\theta}$.
Therefore $g_\theta \in (\lclassg {\infty} (\lsclass {\infty}), X (\lsclass {\infty}))_\theta =
X^\theta (\lsclass {\infty})$ by Proposition~\ref {infint} with
\begin {equation}
\label {gthe}
\|g_\theta\|_{X^\theta (\lsclass {\infty})} \leqslant 2 c_1 \|a\|_{(\BMO, X)_\theta}.
\end {equation}
Observe that by Proposition~\ref {a1interp} applied to $\apclass {1}$-regular lattices $X$ and $X'$
lattices $X^\theta$ and $(X^\theta)' = {X'}^\theta \lclassg {1}^{1 - \theta}$ are also $\apclass {1}$-regular,
so by Proposition~\ref {fsharpest} we have the estimate
\begin {equation}
\label {athe}
\|a\|_{X^\theta} \leqslant c \|a^\sharp\|_{X^\theta}
\end {equation}
for all $a \in X^\theta$ with some $c$ independent of $a$.
Since the function under the supremum in the definition of the Fefferman-Stein maximal function
depends continuously on the coordinates of the vertices of the cube $Q$,
the maximal function $f_\theta^\sharp$ takes the same values if we only take cubes with rational coordinates of the vertices.
Therefore
\begin {multline*}
a^\sharp (x) = f_\theta^\sharp (x) = \sup_j \frac 1 {|Q_j|} \int_{x + Q_j} \left| f_\theta - \frac 1 {|Q|} \int_{x + Q_j} f_\theta \right| =
\sup_j g_{\theta, j} (x)
\end {multline*}
for all $x \in \mathbb R^n$, which links \eqref {gthe} and \eqref {athe} together:
$$
\|a\|_{X^\theta} = \|f_\theta\|_{X^\theta} \leqslant c \|f_\theta^\sharp\|_{X^\theta} = c \|g_\theta\|_{X^\theta (\lsclass {\infty})} \leqslant 2 c_1 \|a\|_{(\BMO, X)_\theta}.
$$
Thus we have verified the claimed continuous inclusion $(\BMO, X)_\theta \subset X^\theta$.  The proof of Theorem~\ref {bintt} is complete.

\subsection* {Acknowledgement}

The author is grateful to S.~V.~Kisliakov who pointed out inaccuracies and provided useful remarks to an early version of this paper.

\bibliographystyle {plain}

\bibliography {bmora}

\end {document}